\documentclass{elsart}
\usepackage{graphics}
\usepackage{graphicx}
\usepackage{epsfig}
\usepackage{amssymb}

\begin{document}
\begin{frontmatter}
\title{Methods for determination and approximation of the domain of attraction}

\author[UVTmat,Paris]{E. Kaslik}
\ead{kaslik@math.univ-paris13.fr}
\author[UVTfiz]{A.M. Balint}
\ead{balint@physics.uvt.ro}
\author[UVTmat]{St. Balint\corauthref{cor}}
\corauth[cor]{Corresponding author.}
\ead{balint@balint.math.uvt.ro}

\address[UVTmat]{Department of Mathematics, West University of
Timi\c{s}oara\\Bd. V. Parvan nr. 4, 300223, Timi\c{s}oara,
Romania\\ phone, fax: +40-256-494002}
\address[UVTfiz]{Department of Physics, West University of
Timi\c{s}oara\\Bd. V. Parvan nr. 4, 300223, Timi\c{s}oara,
Romania}
\address[Paris]{L.A.G.A, UMR 7539, Institut Galil\'{e}e, Universit\'{e} Paris 13\\ 99 Avenue J.B. Cl\'{e}ment,
93430, Villetaneuse, France}

\begin{abstract}
In this paper, an $\mathbb{R}$-analytical function and the
sequence of its Taylor polynomials (which are Lyapunov functions
different from those of Vanelli \& Vidyasagar (1985, Automatica,
21(1):6 9--80)) is presented, in order to determine and
approximate the domain of attraction of the exponentially
asymptotically stable zero steady state of an autonomous,
$\mathbb{R}$-analytical system of differential equations. The
analytical function and the sequence of its Taylor polynomials are
constructed by recurrence formulae using the coefficients of the
power series expansion of $f$ at $0$.
\end{abstract}


\end{frontmatter}

\section{Introduction}
Let be the following system of differential equations:
\begin{equation}\label{dyn.sys}
    \dot{x}=f(x)
\end{equation}
where $f:\mathbb{R}^{n}\rightarrow\mathbb{R}^{n}$ is a function of
class $C^{1}$ on $\mathbb{R}^{n}$ with $f(0)=0$ (i.e. $x=0$ is a
steady state of (\ref{dyn.sys})). If the steady state $x=0$ is
asymptotically stable \cite{Gruyitch}, then the set $D_{a}(0)$ of
all initial states $x^{0}$ for which the solution $x(t;0,x^{0})$
of the initial value problem:
\begin{equation}
    \dot{x}=f(x) \qquad x(0)=x^{0}
\end{equation}
tends to $0$ as $t$ tends to $\infty$, is open and connected and
it is called the domain of attraction (domain of asymptotic
stability \cite{Gruyitch}) of $0$.

The results of Barbashin \cite{Barbashin}, Barbashin-Krasovskii
\cite{Barbashin-Krasovskii} and of Zubov (\cite{Zubov2}, Theorem
19, pp. 52-53, \cite{Zubov}), have probably been the first results
concerning the exact determination of $D_{a}(0)$. In our context,
the theorem of Zubov is the following:

\begin{thm}
An invariant and open set $S$ containing the origin and included
in the hypersphere $B(r)=\{x\in\mathbb{R}^{n}:\|x\|<r\}$, $r>0$,
coincides with the domain of attraction $D_{a}(0)$ if and only if
there exist two functions $V$ and $\psi$ with the following
properties:
\begin{enumerate}
    \item the function $V$ is defined and continuous on $S$, and
    the function $\psi$ is defined and continuous on $\mathbb{R}^{n}$
    \item $-1<V(x)<0$ for any $x\in S\setminus\{0\}$ and $\psi
    (x)>0$, for any $x\in\mathbb{R}^{n}\setminus\{0\}$
    \item $\lim\limits_{x\rightarrow 0}V(x)=0$ and $\lim\limits_{x\rightarrow 0}\psi(x)=0$
    \item for any $\gamma_{2}>0$ small enough, there exist
    $\gamma_{1}>0$ and $\alpha_{1}>0$ such that $V(x)<-\gamma_{1}$ and
    $\psi (x)>\alpha_{1}$, for $\|x\|\geq\gamma_{2}$
    \item for any $y\in\partial S$, $\lim\limits_{x\rightarrow y}V(x)=-1$
    \item $\frac{d}{dt}V(x(t;0;x^{0}))=\psi(x(t;0;x^{0}))[1+V(x(t;0;x^{0}))]$
\end{enumerate}
\label{thm.Zubov}
\end{thm}

\begin{rem}
At this level of generality, the effective determination of
$D_{a}(0)$ using the functions $V$ and $\psi$ from Zubov's theorem
is not possible, because the function $V$ (if $\psi$ is chosen) is
constructed by the method of characteristics, using the solutions
of system (\ref{dyn.sys}). This fact implicitly requests the
knowledge of the domain of attraction $D_{a}(0)$ itself.
\label{rem.Zubov.nefolosibila}
\end{rem}

Another interesting result concerning the exact determination of
$D_{a}(0)$, under the hypothesis that the real parts of the
eigenvalues of the matrix $\frac{\partial f}{\partial x}(0)$ are
negative, is due to Knobloch and Kappel \cite{Knobloch-Kappel}. In
our context, Knobloch-Kappel's theorem is the following:

\begin{thm}
If the real parts of the eigenvalues of the matrix $\frac{\partial
f}{\partial x}(0)$ are negative, then for any function
$\zeta:\mathbb{R}^{n}\rightarrow\mathbb{R}$, with the following
properties:
\begin{enumerate}
    \item $\zeta$ is of class $C^{2}$ on $\mathbb{R}^{n}$
    \item $\zeta(0)=0$ and $\zeta(x)>0$, for any $x\neq 0$
    \item the function $\zeta$ has a positive lower limit on
    every subset of the set $\{x:\|x\|\geq \varepsilon\}$, $\varepsilon >0$
\end{enumerate}
there exists a unique function $V$ of class $C^{1}$ on $D_{a}(0)$
which satisfies
\begin{itemize}
    \item[a.] $\langle\nabla V(x),f(x)\rangle=-\zeta(x)$
    \item[b.] $V(0)=0$
\end{itemize}
In addition, $V$ satisfies the following conditions:
\begin{itemize}
    \item[c.] $V(x)>0$, for any $x\neq 0$
    \item[d.] $\lim\limits_{x\rightarrow y}V(x)=\infty$, for any $y\in\partial
    D_{a}(0)$ or for $\|x\|\rightarrow\infty$
\end{itemize}
\label{thm.Knobloch-Kappel}
\end{thm}

\begin{rem}
The effective determination of $D_{a}(0)$ using the functions $V$
and $\zeta$ from Knobloch-Kappel's theorem (at this level of
generality) is not possible, because the function $V$ (if $\zeta$
is chosen) is constructed by the method of characteristics using
the solutions of system (\ref{dyn.sys}). This fact implicitly
requests the knowledge of $D_{a}(0)$. \label{rem.KK.nefolosibila}
\end{rem}

Vanelli and Vidyasagar have established in
\cite{Vanelli-Vidyasagar} a result concerning the existence of a
maximal Lyapunov function (which characterizes $D_{a}(0)$), and of
a sequence of Lyapunov functions which can be used for
approximating the domain of attraction $D_{a}(0)$. In the context
of our paper, the theorem of Vanelli-Vidyasagar is the following:

\begin{thm}
An open set $S$ which contains the origin coincides with the
domain of asymptotic stability of the asymptotically stable steady
state $x=0$, if and only if there exists a continuous function
$V:S\rightarrow\mathbb{R}_{+}$ and a positive definite function
$\psi$ on $S$ with the following properties:
\begin{enumerate}
    \item $V(0)=0$ and $V(x)>0$, for any $x\in S\setminus\{0\}$ ($V$ is positive definite on $S$)
    \item $D_{r}V(x^{0})=\lim\limits_{t\rightarrow
    0_{+}}\frac{V(x(t;0,x^{0}))-V(x^{0})}{t}=-\psi(x^{0})$, for
    any $x^{0}\in S$
    \item $\lim\limits_{x\rightarrow y}V(x)=\infty$, for any $y\in\partial
    S$ or for $\|x\|\rightarrow\infty$
\end{enumerate}
\label{thm.Vanelli-Vidyasagar}
\end{thm}

\begin{rem}
The determination of $D_{a}(0)$ using the functions $V$ and $\psi$
from Vanelli-Vidyasagar's theorem is not possible, for the same
reason as in the case of the theorems of Zubov and
Knobloch-Kappel. \label{rem.VV.nefolosibila}
\end{rem}

\begin{rem}
Restraining generality, and considering the case of an
$\mathbb{R}$-analytic function $f$, for which the real parts of
the eigenvalues of the matrix $\frac{\partial f}{\partial x}(0)$
are negative, Vanelli and Vidyasagar \cite{Vanelli-Vidyasagar}
establish a second theorem which provides a sequence of Lyapunov
functions, which are not necessarily maximal, but can be used in
order to approximate $D_{a}(0)$. These Lyapunov functions are of
the form:
\begin{equation}
    V_{m}(x)=\frac{r_{2}(x)+r_{3}(x)+...+r_{m}(x)}{1+q_{1}(x)+q_{2}(x)+...+q_{m}(x)}
    \qquad m\in\mathbb{N}
\end{equation}
where $r_{i}$ and $q_{i}$ are $i$-th degree homogeneous
polynomials, constructed using the elements of the matrix
$\frac{\partial f}{\partial x}(0)$, of a positively definite
matrix $G$ and the nonlinear terms from the development of $f$.
The algorithm of the construction of $V_{m}$ is relatively
complex, but does not suppose knowledge of the solutions of system
(\ref{dyn.sys}). \label{rem.VV.Lyap}
\end{rem}

Very interesting results concerning the exact determination of the
domains of attraction (asymptotic stability domains) have been
found by Gruyitch between 1985-1995. These results can be found in
\cite{Gruyitch}, chap. 5. In these results, the function $V$ which
characterizes the domain of attraction is constructed by the
method of characteristics, which uses the solutions of system
(\ref{dyn.sys}). Some illustrative examples are exceptions because
for them $V$ is found in a finite form, for some concrete
functions $f$, but without a precise generally applicable rule.

In the same year as Vanelli and Vidyasagar (1985), Balint
\cite{Balint1}, proved the following theorem:

\begin{thm}(see \cite{Balint1} or \cite{KBB})
If the function $f$ is $\mathbb{R}$-analytic and the real parts of
the eigenvalues of the matrix $\frac{\partial f}{\partial x}(0)$
are negative, then the domain of attraction $D_{a}(0)$ of the
asymptotically stable steady state $x=0$ coincides with the
natural domain of analyticity of the $\mathbb{R}$-analytical
function $V$ defined by
\begin{equation}\label{eq.V}
\langle\nabla V(x),f(x)\rangle=-\|x\|^{2}\qquad V(0)=0
\end{equation}
The function $V$ is strictly positive on $D_{a}(0)\setminus\{0\}$
and $\lim\limits_{x\rightarrow y}V(x)=\infty$ for any
$y\in\partial D_{a}(0)$ or for $\|x\|\rightarrow\infty$.
\label{thm.Balint}
\end{thm}

\begin{rem}
In the case when the matrix $\frac{\partial f}{\partial x}(0)$ is
diagonalizable, recurrence formulae have been established in
\cite{Balint2} (see also \cite{KBB}) for the computation of the
coefficients of the power series expansion in $0$ of the function
$V$ defined by (\ref{eq.V}) (called optimal Lyapunov function in
\cite{Balint2}):

Consider $S:\mathbb{C}^{n}\rightarrow\mathbb{C}^{n}$ an
isomorphism which reduces $\frac{\partial f}{\partial x}(0)$ to
the diagonal form $S^{-1}\frac{\partial f}{\partial
x}(0)S=diag(\lambda_{1},\lambda_{2}...\lambda_{n})$. Let be
$g=S^{-1}\circ f\circ S$ and $W=V\circ S$. If the expansion of $W$
at $0$ is
\begin{equation}\label{serieW}
    W(z_{1},z_{2},...,z_{n})=\sum\limits_{m=2}^{\infty}\sum\limits_{|j|=m}B_{j_{1}j_{2}...j_{n}}z_{1}^{j_{1}}z_{2}^{j_{2}}...z_{n}^{j_{n}}
\end{equation}
\noindent and the expansions at $0$ of the scalar components
$g_{i}$ of $g$ are
\begin{equation}\label{serieg}
    g_{i}(z_{1},z_{2},...,z_{n})=\lambda_{i}z_{i}+\sum\limits_{m=2}^{\infty}\sum\limits_{|j|=m}b^{i}_{j_{1}j_{2}...j_{n}}z_{1}^{j_{1}}z_{2}^{j_{2}}...z_{n}^{j_{n}}
\end{equation}
then the coefficients $B_{j_{1}j_{2}...j_{n}}$ of the development
(\ref{serieW}) are given by the following relations:
\begin{equation}
\label{coefB} B_{j_{1}j_{2}...j_{n}}=
\begin{array}{lll}
\left\{\begin{array}{l}
-\frac{1}{2\lambda_{i_{0}}}\sum\limits_{i=1}^{n}s^{2}_{ii_{0}}
\textrm{
if   } |j|=j_{i_{0}}=2\\ \\
-\frac{2}{\lambda_{p}+\lambda_{q}}\sum\limits_{i=1}^{n}s_{ip}s_{iq}
\textrm{  if   } |j|=2 \textrm{ and } j_{p}=j_{q}=1\\ \\
-\frac{1}{\sum\limits_{i=1}^{n}j_{i}\lambda_{i}}\sum\limits_{p=2}^{|j|-1}\sum\limits_{|k|=p,k_{i}\leq
j_{i}}\sum\limits_{i=1}^{n}[(j_{i}-k_{i}+1)\\
b^{i}_{k_{1}k_{2}...k_{n}}B_{j_{1}-k_{1}...j_{i}-k_{i}+1...j_{n}-k_{n}}]
\textrm{ if } |j|\geq3
\end{array}\right.
\end{array}
\end{equation}
Using these recurrence formulae, the optimal Lyapunov functions
$V$ and the domains of attraction $D_{a}(0)$ for some
two-dimensional systems have been found in \cite{Balint2} in a
finite form. \label{rem.Balint.procedeu}
\end{rem}

\begin{exmp}
\emph{
\begin{equation}\label{ex1}
\begin{array}{ll}
\left\{\begin{array}{l}
\dot{x_{1}}=-\lambda x_{1}+\rho_{1}x_{1}^{2}+\rho_{2}x_{1}x_{2}\\
\dot{x_{2}}=-\lambda x_{2}+\rho_{1}x_{1}x_{2}+\rho_{2}x_{2}^{2}
\end{array}\right.
\end{array}\qquad \lambda>0,\rho_{1},\rho_{1}\in\mathbb{R}^{1}
\end{equation}
The Lyapunov function corresponding to the zero asymptotically
stable steady state of this system is
\begin{equation}
    V(x_{1},x_{2})=\frac{x_{1}^{2}+x_{2}^{2}}{\lambda}[\frac{\lambda^{2}}{(\rho_{1}x_{1}+\rho_{2}x_{2})^{2}}
    \ln\frac{\lambda}{\lambda-(\rho_{1}x_{1}+\rho_{2}x_{2})}-\frac{\lambda}{\rho_{1}x_{1}+\rho_{2}x_{2}}]
\end{equation}
and the domain of attraction is
\begin{equation}
    D_{a}(0)=\{x\in\mathbb{R}^{2}:\rho_{1}x_{1}+\rho_{2}x_{2}<\lambda\}
\end{equation}}
\label{exmp.Balint1}
\end{exmp}

\begin{exmp}
\emph{
\begin{equation}\label{ex1}
\begin{array}{ll}
\left\{\begin{array}{l}
\dot{x_{1}}=-\lambda x_{1}+\rho x_{1}^{3}+\rho x_{1}x_{2}^{2}\\
\dot{x_{2}}=-\lambda x_{2}+\rho x_{1}^{2}x_{2}+\rho x_{2}^{3}
\end{array}\right.
\end{array}\qquad \lambda>0,\rho\in\mathbb{R}^{1}
\end{equation}
The Lyapunov function corresponding to the zero asymptotically
stable steady state of this system is
\begin{equation}
    V(x_{1},x_{2})=\frac{1}{2\rho}\ln\frac{\lambda}{\lambda-\rho(x_{1}^{2}+x_{2}^{2})}
\end{equation}
and the domain of attraction is
\begin{equation}
    D_{a}(0)=\{x\in\mathbb{R}^{2}:\lambda-\rho(x_{1}^{2}+x_{2}^{2})>0\}
\end{equation}}
\label{exmp.Balint2}
\end{exmp}

Therefore, when the function $f$ is $\mathbb{R}$-analytic, the
real parts of the eigenvalues of the matrix $\frac{\partial
f}{\partial x}(0)$ are negative, and the matrix $\frac{\partial
f}{\partial x}(0)$ is diagonalizable, then the optimal Lyapunov
function $V$  can be found theoretically by computing the
coefficients of its power series expansion at $0$, without knowing
the solutions of system (\ref{dyn.sys}). More precisely, in this
way, the "embryo" $V_{0}$ (i.e. the sum of the series) of the
function $V$ is found theoretically on the domain of convergence
$D_{0}$ of the power series expansion. A formula for determining
the region of convergence $D_{0}\subset D_{a}(0)$ of the series of
$V$ can be found in \cite{Balint3} or \cite{KBB}. If $D_{0}$ is a
strict part of $D_{a}(0)$, then the "embryo" $V_{0}$  can be
prolonged using the algorithm of prolongation of analytic
functions:

If $D_{0}$ is strictly contained in $D_{a}(0)$, then there exists
a point $x^{0}\in\partial D_{0}$ such that the function $V_{0}$ is
bounded on a neighborhood of $x^{0}$. Let be a point $x^{0}_{1}\in
D_{0}$ close to $x^{0}$, and the power series development of
$V_{0}$ in $x^{0}_{1}$ (the coefficients of this development are
determined by the derivatives of $V_{0}$ in $x_{1}^{0}$). Using
the formula from \cite{Balint3} or \cite{KBB}, the domain of
convergence $D_{1}$ of the series centered in $x^{0}_{1}$ is
obtained, which gives a new part $D_{1}\setminus(D_{0}\bigcap
D_{1})$ of the domain of attraction $D_{a}(0)$. The sum $V_{1}$ of
the series centered in $x^{0}_{1}$ is a prolongation of the
function $V_{0}$ to $D_{1}$ and coincides with $V$ on $D_{1}$. At
this step, the part $D_{0}\bigcup D_{1}$ of $D_{a}(0)$ and the
restriction of $V$ to $D_{0}\bigcup D_{1}$  are obtained.

If there exists a point $x^{1}\in\partial (D_{0}\bigcup D_{1})$
such that the function $V|_{D_{0}\bigcup D_{1}}$ is bounded on a
neighborhood of $x^{1}$, then the domain $D_{0}\bigcup D_{1}$ is
strictly included in the domain of attraction $D_{a}(0)$. In this
case, the procedure described above is repeated, in a point
$x_{1}^{1}$  close to $x^{1}$.

The procedure cannot be continued in the case when it is found
that on the boundary of the domain $D_{0}\bigcup D_{1}\bigcup
...\bigcup D_{p}$ obtained at step $p$, there are no points having
neighborhoods on which $V|_{D_{0}\bigcup D_{1}\bigcup ...\bigcup
D_{p}}$ is bounded. We illustrate this process in the following
example:

\begin{exmp}
\emph{Consider the following differential equation:
\begin{equation}
    \dot{x}=x(x-1)(x+2)
\end{equation}
$x=0$ is an asymptotically stable steady state for this equation.
The coefficients of the power series development in $0$ of the
optimal Lyapunov function are computed using (\ref{coefB}):
$A_{n}=\frac{2^{n-1}+(-1)^{n}}{3n 2^{n-1}}$, $n\geq 2$. The domain
of convergence $D_{0}=(-1,1)$ of the series is found using the
formula:
\begin{equation}
    x\in D_{0}\qquad \textrm{iff}\qquad \overline{\lim_{n}}\sqrt[n]{|A_{n}x^{n}|}<1
\end{equation}
The embryo $V_{0}(x)$ is unbounded in $1$ and bounded in $-1$, as
$V_{0}(-1)=\frac{\ln 2}{3}$. We expand $V_{0}(x)$ in $-0.9$ close
to $-1$. The coefficients of the series centered in $-0.9$ are:
$A'_{n}=\frac{1}{3n}[\frac{1}{(1.9)^{n}}+\frac{2(-1)^{n}}{(1.1)^{n}}]$.
The domain of convergence $D_{1}$ of the series centered in $-0.9$
is given by:
\begin{equation}
    x\in D_{1}\qquad \textrm{iff}\qquad \overline{\lim_{n}}\sqrt[n]{|A'_{n}(x+0.9)^{n}|}<1
\end{equation}
and it is $D_{1}=(-2,0.2)$. So far, we have obtained the part
$D=D_{0}\bigcup D_{1}=(-2,1)$ of the domain of attraction
$D_{a}(0)$. As the function $V$ is unbounded at both ends of the
interval, we conclude that $D_{a}(0)=(-2,1)$.}
\end{exmp}

We have illustrated how this approximation technique described in
\cite{Balint1,Balint2,Balint3} works in some particular cases. In
more complex cases (for example if the right hand side terms in
(\ref{dyn.sys}) are just polynomials of second degree), we can
only compute effectively the coefficients $A_{j_{1}j_{2}...j_{n}}$
of the expansion of $V$ up to a finite degree $p$. With these
coefficients, the Taylor polynomial of degree $p$ corresponding to
$V$:
\begin{equation}
V_{0}^{p}(x_{1},x_{2},...,x_{n})=\sum\limits_{m=2}^{p}\sum\limits_{|j|=m}A_{j_{1}j_{2}...j_{n}}
x_{1}^{j_{1}}x_{2}^{j_{2}}...x_{n}^{j_{n}}
\end{equation}
can be constructed. In the followings, it will be shown how
$V_{0}^{p}$ can be used in order to approximate $D_{a}(0)$.

\section{Theoretical results}
For $r>0$, we denote by $B(r)=\{x\in\mathbb{R}^{n}:\|x\|<r\}$ the
hypersphere of radius $r$.

\begin{thm}
For any $p\geq 2$, there exists $r_{p}>0$ such that for any $x\in
\overline{B(r_{p})}\setminus\{0\}$ one has:
\begin{enumerate}
    \item $V_{p}(x)>0$
    \item $\langle\nabla V_{p}(x), f(x)\rangle <0$
\end{enumerate}
\label{thm.Vp.lyap.bila}
\end{thm}

\begin{pf}
First, we will prove that for $p=2$, the function $V_{2}$ has the
properties 1. and 2. For this, write the function $f$ as:
\begin{equation}
    f(x)=Ax+g(x)\qquad \textrm{with } A=\frac{\partial f}{\partial x}(0)
\end{equation}
and the equation
\begin{equation}
    \langle\nabla V(x),f(x)\rangle=-\|x\|^{2}
\end{equation}
as
\begin{equation}
     \langle\nabla V_{2}(x),Ax\rangle + \langle\nabla(V-V_{2})(x),Ax+g(x)\rangle + \langle\nabla V_{2}(x),g(x)\rangle =-\|x\|^{2}
\end{equation}
Equating the terms of second degree, we obtain:
\begin{equation}
    \langle\nabla V_{2}(x),Ax\rangle =-\|x\|^{2}
\end{equation}
As $V_{2}(0)=0$, it results that:
\begin{equation}
\label{V2.expresie}
    V_{2}(x)=\int_{0}^{\infty}\|e^{At}x\|^{2}dt
\end{equation}
This shows that $V_{2}(x)>0$ for any
$x\in\mathbb{R}^{n}\setminus\{0\}$.

On the other hand, one has:
\begin{eqnarray}
  \nonumber \langle\nabla V_{2}(x),f(x)\rangle &=& \langle\nabla V_{2}(x),Ax\rangle +\langle\nabla V_{2}(x),g(x)\rangle = \\
  \nonumber &=& -\|x\|^{2} + \langle\nabla V_{2}(x),g(x)\rangle = \\
   &=& -\|x\|^{2}[1-\frac{\langle\nabla V_{2}(x),g(x)\rangle}{\|x\|^{2}}]
\end{eqnarray}
As $\lim\limits_{\|x\|\rightarrow 0}\frac{\langle\nabla
V_{2}(x),g(x)\rangle}{\|x\|^{2}}=0$, there exists $r_{2}>0$ such
that for any $x\in \overline{B(r_{2})}\setminus\{0\}$, we have
$|\frac{\langle\nabla
V_{2}(x),g(x)\rangle}{\|x\|^{2}}|<\frac{1}{2}$. Therefore, for any
$x\in \overline{B(r_{2})}\setminus\{0\}$, we get that:
\begin{equation}
    \langle\nabla V_{2}(x),f(x)\rangle\leq -\frac{1}{2}\|x\|^{2}
\end{equation}
We will show that for any $p>2$, the function $V_{p}$ satisfies
conditions 1. and 2. Write the function $V_{p}$ as
\begin{equation}
    V_{p}(x)=V_{2}(x)[1+\frac{V_{p}(x)-V_{2}(x)}{V_{2}(x)}]\qquad
    x\neq 0
\end{equation}
As $\lim\limits_{\|x\|\rightarrow
0}\frac{V_{p}(x)-V_{2}(x)}{V_{2}(x)}=0$, there exists $r_{p}^{1}$
such that for any $x\in \overline{B(r_{p}^{1})}\setminus\{0\}$, we
have $|\frac{V_{p}(x)-V_{2}(x)}{V_{2}(x)}|<\frac{1}{2}$.
Therefore, for any $x\in \overline{B(r_{p}^{1})}\setminus\{0\}$,
we have:
\begin{equation}
    V_{p}(x)\geq\frac{1}{2}V_{2}(x)>0
\end{equation}
thus, $V_{p}$ satisfies condition 1.

On the other hand, we have:
\begin{eqnarray}
\nonumber \langle\nabla V_{p}(x),f(x)\rangle &=& \langle\nabla
    V_{2}(x),Ax\rangle[1+\frac{\langle\nabla
    (V_{p}-V_{2})(x),f(x)\rangle+\langle\nabla
    V_{2}(x),g(x)\rangle}{\langle\nabla
    V_{2}(x),Ax\rangle}]= \\
   &=& -\|x\|^{2}[1-\frac{\langle\nabla
    (V_{p}-V_{2})(x),f(x)\rangle+\langle\nabla
    V_{2}(x),g(x)\rangle}{\|x\|^{2}}]
\end{eqnarray}
As $\lim\limits_{\|x\|\rightarrow 0}\frac{\langle\nabla
(V_{p}-V_{2})(x),f(x)\rangle+\langle\nabla
V_{2}(x),g(x)\rangle}{\|x\|^{2}}=0$, there exists $r_{p}^{2}$ such
that for any $x\in \overline{B(r_{p}^{2})}\setminus\{0\}$, we have
$|\frac{\langle\nabla (V_{p}-V_{2})(x),f(x)\rangle+\langle\nabla
V_{2}(x),g(x)\rangle}{\|x\|^{2}}|<\frac{1}{2}$. Therefore, for any
$x\in \overline{B(r_{p}^{2})}\setminus\{0\}$, we have:
\begin{equation}
    \langle\nabla V_{p}(x),f(x)\rangle\leq -\frac{1}{2}\|x\|^{2}
\end{equation}
Therefore, for any $x\in \overline{B(r_{p})}\setminus\{0\}$, where
$r_{p}=\min\{r_{p}^{1},r_{p}^{2}\}$, the function $V_{p}$
satisfies conditions 1. and 2.
\end{pf}

\begin{cor}
For any $p\geq 2$, there exists a maximal domain
$G_{p}\subset\mathbb{R}^{n}$ such that $0\in G_{p}$ and for any
$x\in G_{p}\setminus\{0\}$, function $V_{p}$ verifies 1. and 2.
from Theorem \ref{thm.Vp.lyap.bila}. In other words, for any
$p\geq 2$ the function $V_{p}$ is a Lyapunov function for
(\ref{dyn.sys}) (in the sense of \cite{Gruyitch}).
\end{cor}

\begin{rem}
Theorem \ref{thm.Vp.lyap.bila} provides that the Taylor
polynomials of degree $p\geq 2$ associated to $V$ in $0$ are
Lyapunov functions. This sequence of Lyapunov functions is
different of that provided by Vanelli and Vidyasagar in
\cite{Vanelli-Vidyasagar}.
\end{rem}

\begin{thm}
For any $p\geq 2$, there exists $c>0$ and a closed and connected
set $S$ of points from $x\in\mathbb{R}^{n}$, with the following
properties:
\begin{enumerate}
    \item $0\in Int(S)$
    \item $V_{p}(x)< c$ for any $x\in Int(S)$
    \item $V_{p}(x)=c$ for any $x\in\partial S$
    \item $S$ is compact and included in the set $G_{p}$.
\end{enumerate}
\label{thm.Ncp.exist}
\end{thm}

\begin{pf}
Let be $p\geq 2$ and $r_{p}>0$ determined in Theorem
\ref{thm.Vp.lyap.bila}. Let be
$c=\min\limits_{\|x\|=r_{p}}V_{p}(x)$ and
$S'=\{x\in\overline{B(r_{p})}:V_{p}(x)< c\}$. It is obvious that
$c>0$ and that there exist $x^{\star}$ with $\|x^{\star}\|=r_{p}$
such that $V(x^{\star})=c$. The set $S'$ is open, $0\in S'$ and
$S'\subset \overline{B(r_{p})}\subset G_{p}$.

We will prove that $V_{p}(x)=c$ for any $x\in\partial S'$. Let be
$\bar{x}\in\partial S'$. Thus, $\|\bar{x}\|\leq r_{p}$ and there
exists a sequence $x^{k}\in S'$ such that $x^{k}\rightarrow
\bar{x}$ as $k\rightarrow\infty$. As $V_{p}(x^{k})<c$, we have
that
$V_{p}(\bar{x})=\lim\limits_{k\rightarrow\infty}V_{p}(x^{k})\leq
c$. The case $\|\bar{x}\|= r_{p}$ and $V_{p}(\bar{x})<c$ is
impossible, because $c=\min\limits_{\|x\|=r_{p}}V_{p}(x)$. The
case $\|\bar{x}\|< r_{p}$ and $V_{p}(\bar{x})<c$ is also
impossible, because this would mean that $\bar{x}$ belongs to the
interior of the set $S'$, and not to its boundary. Therefore, for
any $\bar{x}\in\partial S'$ we have $V_{p}(\bar{x})=c$.

If the set $S'$ is not connected (see Example
\ref{ex.nu.rad.cresc} in this paper), we denote by $S''$ its
connected component which contains the origin, and let be
$S=\overline{S''}$. Then it is obvious that $S$ is connected
(being the closure of the open connected set $S''$), $0\in
Int(S)=S''$, and that for any $x\in Int(S)=S''$, we have
$V_{p}(x)<c$. More, as $\partial S=\partial S''$, we have
$V_{p}(x)=c$ for any $x\in
\partial S$. As $S''$ is bounded, we obtain that the closed set
$S$ is also bounded, thus, it is compact. As $S''\subset
\overline{B(r_{p})}\subset G_{p}$, we have that
$S=\overline{S''}\subset \overline{B(r_{p})}\subset G_{p}$.
Therefore, $S$ satisfies the properties 1-4.
\end{pf}

\begin{lem}
Let be $p\geq 2$, $c>0$ and a closed and connected set $S$
satisfying 1-4 from Theorem \ref{thm.Ncp.exist}. Then for any
$x^{0}\in S$, the solution $x(t;0,x^{0})$ of system
(\ref{dyn.sys}) starting from $x^{0}$ is defined on $[0,\infty)$
and belongs to $Int(S)$ for any $t>0$. \label{lemma.Npc.prelung}
\end{lem}

\begin{pf}
Let be $x^{0}\in S$. We denote by $[0,\beta_{x^{0}})$ the right
maximal interval of existence of the solution $x(t;0,x^{0})$ of
system (\ref{dyn.sys}) with starting state $x^{0}$.

First, if $x^{0}\in Int(S)\setminus\{0\}$, we show that
$x(t;0,x^{0})\in Int(S)$, for all $t\in [0,\beta_{x^{0}})$.
Suppose the contrary, i.e. there exists $T\in(0,\beta_{x^{0}})$
such that $x(t;0,x^{0})\in Int(S)$, for $t\in [0,T)$ and
$x(T;0,x^{0})\in \partial S$ (i.e. $V_{p}(x(T;0,x^{0}))=c$). As
$x(t;0,x^{0})\in G_{p}\setminus\{0\}$, for $t\in [0,T)$,
$V_{p}(x(t;0,x^{0}))$ is strictly decreasing, and it follows that
$V_{p}(x(t;0,x^{0}))<V_{p}(x^{0})<c$, for $t\in (0,T)$. Therefore
$V_{p}(x(T;0,x^{0}))<c$, which contradicts the supposition
$x(T;0,x^{0})\in \partial S$. Thus, $x(t;0,x^{0})\in Int(S)$, for
all $t\in [0,\beta_{x^{0}})$. (It is clear that for $x^{0}=0$, the
solution $x(t;0,0)=0\in Int(S)$, for all $t\geq 0$.)

If $x^{0}\in \partial S$, we show that $x(t;0,x^{0})\in Int(S)$,
for all $t\in (0,\beta_{x^{0}})$. As the compact set $S$ is a
subset of the domain $G_{p}$, the continuity of $x(t;0,x^{0})$
provides the existence of $T_{x^{0}}>0$ such that $x(t;0,x^{0})\in
G_{p}\setminus\{0\}$ for any
$t\in[0,T_{x^{0}}]\subset[0,\beta_{x^{0}})$. Therefore
$V_{p}(x(t;0,x^{0}))$ is strictly decreasing on $[0,T_{x^{0}}]$,
and it follows that $V_{p}(x(t;0,x^{0}))<V_{p}(x^{0})= c$, for any
$t\in (0,T_{x^{0}})$. This means that $V_{p}(x(t;0,x^{0}))\in
Int(S)$, for any $t\in (0,T_{x^{0}}]$. The first part of the proof
guarantees that $x(t;0,x^{0})\in Int(S)$, for all $t\in
[T_{x^{0}},\beta_{x^{0}})$, therefore, for all $t\in
(0,\beta_{x^{0}})$.

In conclusion, for any $x^{0}\in S$, we have that $x(t;0,x^{0})\in
Int(S)$, for all $t\in (0,\beta_{x^{0}})$.

As for any $x^{0}\in S$, the solution $x(t;0,x^{0})$ defined on
$[0,\beta_{x^{0}})$, belongs to the compact $S$, we obtain that
$\beta_{x^{0}}=\infty$ and the solution $x(t;0,x^{0})$ is defined
on $[0,\infty)$, for each $x^{0}\in S$. More, $x(t;0,x^{0})\in
Int(S)$, for all $t>0$.
\end{pf}

\begin{rem}
Lemma \ref{lemma.Npc.prelung} states that a closed and connected
set $S$ satisfying 1-4 from Theorem \ref{thm.Ncp.exist} is
positively invariant to the flow of system (\ref{dyn.sys}).
\end{rem}

\begin{thm}(LaSalle-type theorem)
Let be $p\geq 2$, $c>0$ and a closed and connected set $S$
satisfying 1-4 from Theorem \ref{thm.Ncp.exist}. Then $S$ is a
part of the of the domain of attraction $D_{a}(0)$.
\label{thm.Ncp.parte.DA}
\end{thm}

\begin{pf}
Let be $x^{0}\in S\setminus\{0\}$. To prove that
$\lim\limits_{t\rightarrow\infty}x(t;0,x^{0})=0$, it is sufficient
to prove that
$\lim\limits_{k\rightarrow\infty}x(t_{k};0,x^{0})=0$, for any
sequence $t_{k}\rightarrow\infty$.

Consider $t_{k}\rightarrow\infty$. The terms of the sequence
$x(t_{k};0,x^{0})$ belong to the compact $S$. Thus, there exits a
convergent subsequence $x(t_{k_{j}};0,x^{0})\rightarrow y^{0}\in
S$.

It can be shown that
\begin{equation}\label{ineg0}
    V_{p}(x(t;0,x^{0}))\geq V_{p}(y^{0}) \textrm{ for all } t\geq 0
\end{equation}
For this, observe that $x(t_{k_{j}};0,x^{0})\rightarrow y^{0}$ and
$V_{p}$ is strictly decreasing along the trajectories, which
implies that $V_{p}(x(t_{k_{j}};0,x^{0}))\geq V_{p}(y^{0})$ for
any $k_{j}$. On the other hand, for any $t\geq 0$, there exists
$k_{j}$ such that $t_{k_{j}}\geq t$, and therefore
$V_{p}(x(t;0,x^{0}))\geq V_{p}(x(t_{k_{j}};0,x^{0}))\geq
V_{p}(y^{0})$.

We show now that $y^{0}=0$. Suppose the contrary, i.e. $y^{0}\neq
0$. Inequality (\ref{ineg0}) becomes
\begin{equation}\label{contrad}
    V_{p}(x(t;0,x^{0}))\geq V_{p}(y^{0})>0 \textrm{ for all } t\geq 0
\end{equation}
As $V_{p}(x(s;0,y^{0}))$ is strictly decreasing on $[0,\infty)$,
we find that
\begin{equation}
    V_{p}(x(s;0,y^{0}))<V_{p}(y^{0})\textrm{ for all } s>0
\end{equation}
For $\bar{s}>0$, there exists a neighborhood
$U_{x(\bar{s};0,y^{0})}\subset S$ of $x(\bar{s};0,y^{0})$ such
that for any $x\in U_{x(\bar{s};0,y^{0})}$ we have
$0<V_{p}(x)<V_{p}(y^{0})$. On the other hand, for the neighborhood
$U_{x(\bar{s};0,y^{0})}$ there exists a neighborhood
$U_{y^{0}}\subset S$ of $y^{0}$ such that $x(\bar{s};0,y)\in
U_{x(\bar{s};0,y^{0})}$ for any $y\in U_{y^{0}}$. Therefore:
\begin{equation}\label{ineg1}
    V_{p}(x(\bar{s};0,y))<V_{p}(y^{0})\textrm{ for all } y\in U_{y^{0}}
\end{equation}
As $x(t_{k_{j}};0,x^{0})\rightarrow y^{0}$, there exists
$k_{\bar{j}}$ such that $x(t_{k_{j}};0,x^{0})\in  U_{y^{0}}$, for
any $k_{j}\geq k_{\bar{j}}$. Making $y=x(t_{k_{j}};0,x^{0})$ in
(\ref{ineg1}), it results that
\begin{equation}\label{ineg2}
    V_{p}(x(\bar{s}+t_{k_{\bar{j}}};0,x^{0}))=V_{p}(x(\bar{s};0,x(t_{k_{\bar{j}}};0,x^{0})))<V_{p}(y^{0})
    \qquad \textrm{for }k_{j}\geq k_{\bar{j}}
\end{equation}
which contradicts (\ref{contrad}). This means that $y^{0}=0$,
consequently, every convergent subsequence of $x(t_{k};0,x^{0})$
converges to $0$. This provides that the sequence
$x(t_{k};0,x^{0})$ is convergent to $0$, for any $t_{k}\rightarrow
\infty$, thus $\lim\limits_{t\rightarrow\infty}x(t;0,x^{0})=0$,
and $x^{0}\in D_{a}(0)$.

Therefore, the set $S$ is contained in the domain of attraction of
$D_{a}(0)$.
\end{pf}

\begin{cor}
For any $p\geq 2$ and $c>0$ there exists at most one closed and
connected set satisfying 1-4 from Theorem \ref{thm.Ncp.exist}.
\end{cor}

\begin{pf}
Suppose the contrary, i.e. for a $p\geq 2$ and $c>0$ there exist
two different closed and connected sets $S_{1}$ and $S_{2}$
satisfying 1-4 from Theorem \ref{thm.Ncp.exist}. Assume for
example that there exists $x^{0}\in S_{1}\setminus S_{2}$. Due to
Theorem \ref{thm.Ncp.parte.DA}, $S_{1}\subset D_{a}(0)$ and
therefore $\lim\limits_{t\rightarrow\infty}x(t;0,x^{0})=0$. As
$x^{0}\notin S_{2}$, and $S_{2}$ is a closed and connected
neighborhood of $0$, there exists $T>0$ such that
$x(T;0,x^{0})\in\partial S_{2}$. Therefore,
$V_{p}(x(T;0,x^{0}))=c$ which contradicts Lemma
\ref{lemma.Npc.prelung}. Consequently, we have $S_{1}\subseteq
S_{2}$. By the same reasons, $S_{2}\subseteq S_{1}$. Finally,
$S_{1}= S_{2}$.
\end{pf}

\begin{rem}
If for $p\geq 2$ and $c>0$ there exists a closed and connected set
satisfying 1-4 from Theorem \ref{thm.Ncp.exist}, then it is unique
and it will be denoted by $N_{p}^{c}$.
\end{rem}

\begin{cor}
Any set $N_{p}^{c}$ is included in the domain of attraction
$D_{a}(0)$.
\end{cor}

\begin{lem}
Let be $p\geq 2$ and $c>0$ such that there exists the set
$N_{p}^{c}$. Then, for any $c'\in(0,c]$ the set $\{x\in
N_{p}^{c}:V_{p}(x)\leq c'\}$ coincides with the set $N_{p}^{c'}$.
\label{lemma.Npc.incluz.conex}
\end{lem}

\begin{pf}
Let be $c'\in(0,c]$. It is obvious that $N_{p}^{c'}$ is included
in the set $\{x\in N_{p}^{c}:V_{p}(x)\leq c'\}$. Let be $x^{0}\in
N_{p}^{c}$ such that $V_{p}(x^{0})\leq c'$.  We know that
$V_{p}(x(t;0,x^{0}))<V_{p}(x^{0})\leq c'$, for any $t>0$. Theorem
\ref{thm.Ncp.parte.DA} provides that $x^{0}\in N_{p}^{c}\subset
D_{a}(0)$, therefore, $x^{0}$ is connected to $0$ through the
continuous trajectory $x(t;0,x^{0})$, along which $V_{p}$ takes
values below $c'$. In conclusion, $x^{0}\in N_{p}^{c'}$.
\end{pf}

\begin{thm}
If for $p\geq 2$ and $c>0$ there exists $N_{p}^{c}$, then for any
$c'\in(0,c)$ there exists $N_{p}^{c'}$ and $N_{p}^{c'}\subset
N_{p}^{c}$. More, for any $c_{1},c_{2}\in (0,c)$ we have
$N_{p}^{c_{1}}\subset N_{p}^{c_{2}}$ if and only if
$c_{1}<c_{2}$.\label{thm.Ncp.incluziune}
\end{thm}

\begin{pf}
Lemma \ref{lemma.Npc.incluz.conex} provides that for any
$c'\in(0,c)$ there exists $N_{p}^{c'}=\{x\in
N_{p}^{c}:V_{p}(x)\leq c'\}$. It is obvious that
$N_{p}^{c'}\subset N_{p}^{c}$.

Let's show that for any $c_{1},c_{2}\in (0,c)$ we have
$N_{p}^{c_{1}}\subset N_{p}^{c_{2}}$ if and only if $c_{1}<c_{2}$.

To show the \emph{necessity}, let's suppose the contrary, i.e.
$N_{p}^{c_{1}}\subset N_{p}^{c_{2}}$ and $c_{1}\geq c_{2}$. Let be
$x^{0}\in\partial N_{p}^{c_{2}}\subset G_{p}$. Then
$V_{p}(x^{0})=c_{2}$ and as $x^{0}\in G_{p}$, we get that
\begin{equation}\label{c2.contrad}
    V_{p}(x(t;0,x^{0}))\leq V_{p}(x^{0})=c_{2}\qquad \textrm {for any
    }t\geq 0
\end{equation}
Theorem \ref{thm.Ncp.parte.DA} provides that $x^{0}\in\partial
N_{p}^{c_{2}}\subset D_{a}(0)$, therefore
$\lim\limits_{t\rightarrow\infty}x(t;0,x^{0})=0$. As
$N_{p}^{c_{1}}$ and $N_{p}^{c_{2}}$ are connected neighborhoods of
$0$ and $N_{p}^{c_{1}}\subset N_{p}^{c_{2}}$, there exists $T\geq
0$ such that $x(T;0,x^{0})\in \partial N_{p}^{c_{1}}$. This means
that $V_{p}(x(T;0,x^{0}))=c_{1}\geq c_{2}$, and (\ref{c2.contrad})
provides that $c_{1}=c_{2}$. As $N_{p}^{c_{1}}$ is strictly
included in $N_{p}^{c_{2}}$, there exists $\bar{x}\in\partial
N_{p}^{c_{1}}$ (i.e. $V_{p}(\bar{x})=c_{1}=c_{2}$) such that
$\bar{x}\in Int(N_{p}^{c_{2}})$. This contradicts the property 2
from Theorem \ref{thm.Ncp.exist} concerning $N_{p}^{c_{2}}$. In
conclusion, $c_{1}< c_{2}$.

To prove the \emph{sufficiency}, let's suppose that $c_{1}<c_{2}$
and let be $x^{0}\in N_{p}^{c_{1}}\setminus\{0\}$. As $x^{0}\in
N_{p}^{c_{1}}\subset D_{a}(0)$, we have that
$\lim\limits_{t\rightarrow\infty}x(t;0,x^{0})=0$, so $x^{0}$ is
connected to $0$ through the continuous trajectory $x(t;0,x^{0})$.
More, as $x^{0}\in N_{p}^{c_{1}}\setminus\{0\}$, we have
$V_{p}(x^{0})\leq c_{1}\leq c_{2}$. This means that $x^{0}\in
N_{p}^{c_{2}}$, therefore $N_{p}^{c_{1}}\subseteq N_{p}^{c_{2}}$.
The inclusion is strict, because $N_{p}^{c_{1}}=N_{p}^{c_{2}}$
means $\partial N_{p}^{c_{1}}=\partial N_{p}^{c_{2}}$, i.e.
$c_{1}=c_{2}$, which contradicts $c_{1}<c_{2}$.
\end{pf}

\begin{cor}
For a given $p\geq 2$, the set of all $N_{p}^{c}$-s is totally
ordered and $\bigcup\limits_{c}N_{p}^{c}$ is included in
$D_{a}(0)$. Therefore, for a given $p\geq 2$, the largest part of
$D_{a}(0)$ which can be found by this method is
$\bigcup\limits_{c}N_{p}^{c}$.
\end{cor}

For any $p\geq 2$ let be $R_{p}=\{r>0:\overline{B(r)}\subset
G_{p}\}$. For $r\in R_{p}$ we denote by
$c_{p}^{r}=\inf\limits_{\|x\|=r}V_{p}(x)$.

\begin{cor}
For any $r\in R_{p}$, there exists the set $N_{p}^{c_{p}^{r}}$ and
$N_{p}^{c_{p}^{r}}\subseteq \overline{B(r)}$.
\label{cor.Ncpcpr.exist}
\end{cor}

\begin{cor}
For any $p\geq 2$ and any $r',r''\in R_{p}$, $r'<r''$, such that
$V_{p}$ is radially increasing on $\overline{B(r'')}$ we have
$c_{p}^{r'}<c_{p}^{r''}$. \label{cor.Ncpcpr.incluziune}
\end{cor}

\begin{rem} In some cases, it can be shown that the function
$V_{p}$ is radially increasing on $G_{p}$:
\begin{itemize}
    \item[a.] $V_{2}$ is radially increasing on $\mathbb{R}^{n}$;
    \item[b.] If $n=1$, then for any $p\geq 2$, $V_{p}$ is
    radially increasing on $G_{p}$.
\end{itemize}
This result is not true in general, provided by the following
example:
\end{rem}

\begin{exmp}
\emph{Let be the following system of differential equations:
\begin{equation}
\begin{array}{ll}
\left\{\begin{array}{l}
\dot{x_{1}}=-x_{1}-x_{1}x_{2}\\
\dot{x_{2}}=-x_{2}+x_{1}x_{2}
\end{array}\right.
\end{array}
\end{equation}
for which $(0,0)$ is an asymptotically stable steady state. For
$p=3$ the Lyapunov function $V_{3}(x_{1},x_{2})$ is given by:
\begin{equation}
    V_{3}(x_{1},x_{2})=\frac{1}{2}(x_{1}^{2}+x_{2}^{2})+\frac{1}{3}(x_{1}x_{2}^{2}-x_{2}x_{1}^{2})
\end{equation}
Consider the point $(3\sqrt{5},\sqrt{5})\in \partial G_{3}$ and
let be $g:[0,1)\rightarrow G_{3}$ defined by
$g(\lambda)=V_{3}(3\sqrt{5}\lambda,\sqrt{5}\lambda)$. The function
$g$ is increasing on $[0,\frac{\sqrt{5}}{3}]$ and decreasing on
$(\frac{\sqrt{5}}{3},1)$, therefore, the Lyapunov function $V_{3}$
is not radially increasing on the direction
$(3\sqrt{5},\sqrt{5})$. In conclusion, $V_{3}$ is not radially
increasing on $G_{3}$.}

\emph{More, for $c=0.32$ there exists $N_{3}^{c}$, but the set
$\{x=(x_{1},x_{2})\in G_{3}:V_{3}(x_{1},x_{2})\leq c\}$ is not
connected. The reason is that the point
$(\bar{x_{1}},\bar{x_{2}})=(\frac{123}{8},\frac{41}{24})\in\partial
G_{3}$ with $V_{3}(\bar{x_{1}},\bar{x_{2}})=0$ has a nonempty
neighborhood $U$ such that $V_{3}(x_{1},x_{2})\leq c$, for any
$(x_{1},x_{2})\in G_{3}\bigcap U$ and $(G_{3}\bigcap U)\bigcap
N_{3}^{c}= \emptyset$.} \label{ex.nu.rad.cresc}
\end{exmp}

\begin{thm}
For any $p\geq 2$ there exists $\rho_{p}>0$ such that $V_{p}$ is
radially increasing on $\overline{B(\rho_{p})}$.
\end{thm}

\begin{pf}
It can be easily verified that $V_{2}$ is radially increasing on
$\mathbb{R}^{n}$, using relation (\ref{V2.expresie}). This
provides that for any $x\in\mathbb{R}^{n}\setminus\{0\}$, the
function $g_{2}^{x}:\mathbb{R}_{+}\rightarrow\mathbb{R}_{+}$
defined by $g_{2}^{x}(\lambda)=V_{2}(\lambda x)$ is strictly
increasing on $\mathbb{R}_{+}$, therefore
$\frac{d}{d\lambda}g_{2}^{x}(\lambda)>0$ on
$\mathbb{R}_{+}^{\star}$, i.e.
\begin{equation}\label{V2.ec.deriv}
\langle\nabla V_{2}(\lambda x),x \rangle>0\qquad\textrm{for any }
\lambda>0\textrm{ and }x\in\mathbb{R}^{n}\setminus\{0\}
\end{equation}
Let be $p>2$, $x\in\mathbb{R}^{n}\setminus\{0\}$ and
$g_{p}^{x}:\mathbb{R}_{+}\rightarrow\mathbb{R}_{+}$ defined by
$g_{p}^{x}(\lambda)=V_{p}(\lambda x)$. One has:
\begin{eqnarray}
\nonumber  \frac{d}{d\lambda}g_{p}^{x}(\lambda) &=& \langle\nabla
V_{p}(\lambda x),x \rangle=\langle\nabla V_{2}(\lambda x),x \rangle+\langle\nabla (V_{p}-V_{2})(\lambda x),x \rangle =\\
  &=&\langle\nabla V_{2}(\lambda x),x \rangle(1+\frac{\langle\nabla (V_{p}-V_{2})(\lambda x),x
  \rangle}{\langle\nabla V_{2}(\lambda x),x \rangle})
\label{Vp.ec.deriv}
\end{eqnarray}
As $\lim\limits_{x\rightarrow
0}\frac{\langle\nabla(V_{p}-V_{2})(\lambda
x),x\rangle}{\langle\nabla V_{2}(\lambda x),x \rangle}=0$, there
exists $\rho_{p}>0$ such that $|\frac{\langle\nabla
(V_{p}-V_{2})(\lambda x),x\rangle}{\langle\nabla V_{2}(\lambda
x),x \rangle}|\leq \frac{1}{2}$, for any $x\in
\overline{B(\rho_{p})}\setminus\{0\}$. Relation
(\ref{Vp.ec.deriv}) provides that for any $x\in
\overline{B(\rho_{p})}\setminus\{0\}$, we have:
\begin{equation}
    \frac{d}{d\lambda}g_{p}^{x}(\lambda)\geq\frac{1}{2}\langle\nabla V_{2}(\lambda x),x
    \rangle>0\qquad\textrm{ for any }\lambda>0
\end{equation}
Therefore, for any $x\in \overline{B(\rho_{p})}\setminus\{0\}$,
the function $g_{p}^{x}$ is strictly increasing on
$\mathbb{R}^{n}$, i.e. $V_{p}$ is radially increasing on
$\overline{B(\rho_{p})}$.
\end{pf}

\begin{thm}
Let be $p\geq 2$ and $c>0$ such that there exists the set
$N_{p}^{c}$. Suppose that for any $c'\leq c$, the sets
$N_{p}^{c'}$ have the star-property, i.e. for any $x\in
N_{p}^{c'}$ and for any $\lambda\in [0,1)$ one has $\lambda x\in
Int(N_{p}^{c'})$. Then $V_{p}$ is radially increasing on
$N_{p}^{c}$.
\end{thm}

\begin{pf}
Let be $x^{0}\in \partial N_{p}^{c}$ and
$0<\lambda_{1}<\lambda_{2}\leq 1$. We have to show that
$V_{p}(\lambda_{1}x^{0})<V_{p}(\lambda_{2}x^{0})$. Denote
$c_{1}=V_{p}(\lambda_{1}x^{0})>0$,
$c_{2}=V_{p}(\lambda_{2}x^{0})>0$ and suppose the contrary, i.e.
$c_{1}\geq c_{2}$. Theorem \ref{thm.Ncp.incluziune} provides that
$N_{p}^{c_{1}}\supseteq N_{p}^{c_{2}}$. Lemma
\ref{lemma.Npc.incluz.conex} guarantees that
$\lambda_{2}x\in\partial N_{p}^{c_{2}}$. As $N_{p}^{c_{2}}$ has
the star-property, then for
$\lambda=\frac{\lambda_{1}}{\lambda_{2}}\in (0,1)$, we have that
$\lambda (\lambda_{2}x)=\lambda_{1}x\in Int(N_{p}^{c_{2}})$, so
$c_{1}=V_{p}(\lambda_{1}x)< c_{2}$ which contradicts the
supposition $c_{1}\geq c_{2}$. Therefore, $V_{p}$ is radially
increasing on $N_{p}^{c}$.
\end{pf}

\begin{rem}
\begin{itemize}
    \item[a.] For any $x\in D_{0}$, there exists $p_{x}\geq 2$ such that $x\in
G_{p}$ for any $p\geq p_{x}$;
    \item[b.] If $n=1$, there exists $p_{0}\geq 2$ such that $D_{0}\subset
    G_{p}$, for any $p\geq p_{0}$.
    \item[c.] If there exists $r>0$ such that $\overline{B(r)}\subset
    G_{p}$ for any $p\geq 2$, then there exists $p_{0}\geq 2$ such that $D_{0}\subset
    G_{p_{0}}$.
\end{itemize}
\end{rem}

\begin{conj}
For any $x\in D_{a}(0)$ there exists $p\geq 2$ and $c>0$ such that
$x\in N_{p}^{c}$.
\end{conj}

\section{Numerical example: the Van der Pol system}

We consider the following system of differential equations:

\begin{equation}\label{sys.Van.der.Pol}
\begin{array}{ll}
\left\{\begin{array}{l}
\dot{x_{1}}=-x_{2}\\
\dot{x_{2}}=x_{1}-x_{2}+x_{1}^{2}x_{2}
\end{array}\right.
\end{array}
\end{equation}

The $(0,0)$ steady state of (\ref{sys.Van.der.Pol}) is
asymptotically stable. The boundary of the domain of attraction of
$(0,0)$ is a limit cycle of (\ref{sys.Van.der.Pol}).

For $p=20$ we have computed that the largest value $c>0$ for which
there exists the set $N_{p}^{c}$ is $c_{20}=8.8466$. For $p=50$,
the largest value $c>0$ for which there exists the set $N_{p}^{c}$
is $c_{50}=13.887$. In the figures below, the thick black curve
represents the boundary of $D_{a}(0,0)$, the thin black curve
represents the boundary of $G_{p}$ and the gray surface represents
the set $N_{p}^{c_{p}}$. The set $N_{p}^{c_{50}}$ approximates
very well the domain of attraction of $(0,0)$.

\begin{figure}[htbp]
\begin{minipage}[t]{0.5\linewidth}
\centering
\includegraphics*[bb=3cm 0cm 13.5cm
10.5cm,width=7cm,angle=0]{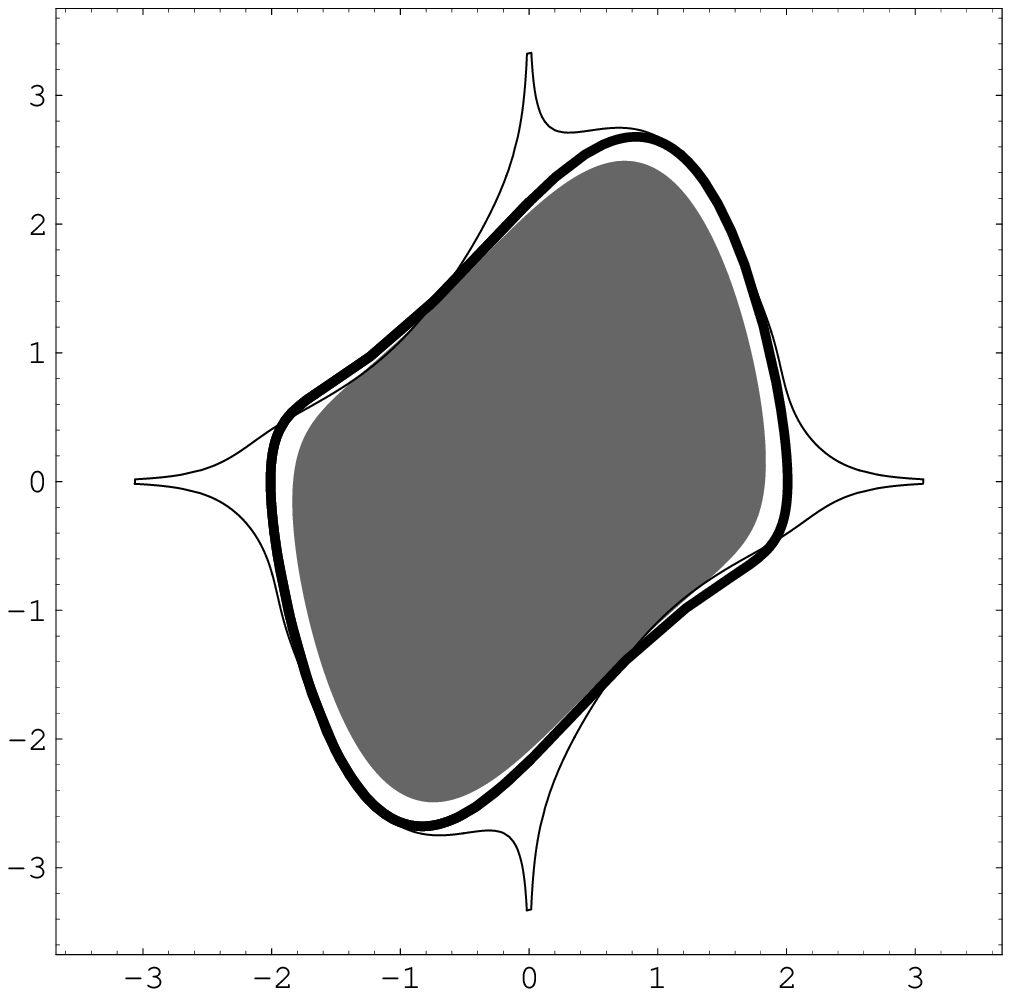}\caption{The sets
$N_{p}^{c_{20}}$, $G_{20}$ and $D_{a}(0,0)$ for system
(\ref{sys.Van.der.Pol})}
\end{minipage}
\begin{minipage}[t]{0.5\linewidth} \centering
\includegraphics*[bb=3cm 0cm 13.5cm
10.5cm,width=7cm,angle=0]{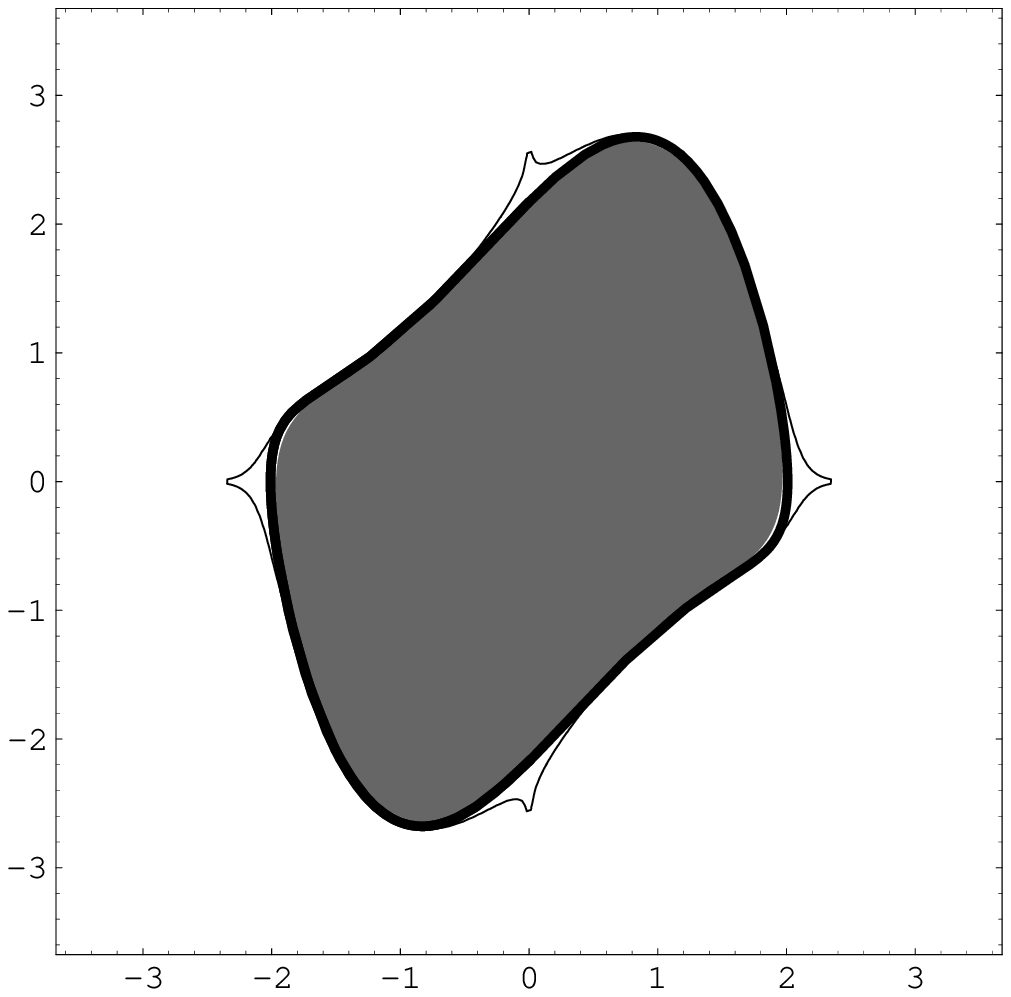}\caption{The sets
$N_{p}^{c_{50}}$, $G_{50}$ and $D_{a}(0,0)$ for system
(\ref{sys.Van.der.Pol})}
\end{minipage}
\end{figure}

\end{document}